\documentclass[12pt]{article}
\usepackage{amsmath}
\usepackage{amsfonts}
\usepackage{amsthm}
\usepackage[dutch,english]{babel}

\newtheorem{theorem}{Theorem}[section]
\newtheorem{lemma}[theorem]{Lemma}
\newtheorem{corollary}[theorem]{Corollary}
\newtheorem{remark}[theorem]{Remark}
\def \qbinom#1#2{\genfrac{[}{]}{0pt}{}{#1}{#2}}
\numberwithin{equation}{section}

\begin{document}
\title{Irrationality proof of a $q$-extension of $\zeta(2)$ \\ using little
        $q$-Jacobi polynomials}
\author{Christophe Smet and Walter Van Assche \\
    \em Department of Mathematics, Katholieke Universiteit Leuven, Belgium}
\date{\today}
\maketitle
\begin{abstract}
We show how one can use Hermite-Pad\'{e} approximation and little
$q$-Jacobi polynomials to construct rational approximants for
$\zeta_q(2)$.  These numbers are $q$-analogues of the
well known $\zeta(2)$.  Here $q=\frac{1}{p}$, with $p$ an
integer greater than one.  These approximants are good enough to
show the irrationality of $\zeta_q(2)$ and they allow
us to calculate an upper bound for its measure of irrationality:
$\mu\left(\zeta_q(2)\right)\leq 10\pi^2/(5\pi^2-24) \approx 3.8936$.  This is sharper than
the upper bound given by Zudilin (\textit{On the irrationality measure for
a $q$-analogue of $\zeta(2)$}, Mat.\ Sb.\  \textbf{193} (2002), no. 8, 49--70).
\end{abstract}

Key words: little $q$-Jacobi polynomials, irrationality, measure
of irrationality, $q$-zeta function

\section{Introduction}

The $\zeta$-function at integer points,
$\zeta(s)=\sum_{k=1}^\infty 1/k^s$ $(s \in \mathbb{N})$ has a
$q$-analogue, defined by
\[\zeta_q(s)=\sum_{k=1}^\infty
\frac{k^{s-1}q^k}{1-q^k}.\]
It makes sense to call this a
$q$-analogue, since
\[\lim_{q\rightarrow1^-}(1-q)^s\zeta_q(s)=(s-1)!\zeta(s).\]
One property this $\zeta_q(s)$ shares with $\zeta(s)$ is that a
lot of questions concerning irrationality remain to be answered.
So far, for $q=1/p$ with $p\in \mathbb{N}\setminus\{0,1\}$, only
$\zeta_q(1)$ and $\zeta_{q}(2)$ have been shown to be irrational.
The former was done by Borwein \cite{bor1a},\cite{bor1b} in 1991,
and, using a different approach, by Bundschuh and
V\"{a}\"{a}n\"{a}nen \cite{buva} in 1994.  Note that a result by
B\'{e}zivin \cite{bezivin} of the year 1988 can be used to prove
this irrationality.  The irrationality of $\zeta_q(2)$ was proven
by Duverney \cite{duv} in 1995, the transcendence of $\zeta_q(2)$
(and in fact of $\zeta_q(2s), s\in\mathbb{N}$) is a consequence of
a general result by Nesterenko \cite{nes}.  Moreover, the three
values $1,\zeta_q(1),\zeta_{q}(2)$ have been shown to be
$\mathbb{Q}$-linearly independent by Postelmans and Van Assche
\cite{kelly}. In this paper we will prove the following result.

\begin{theorem}\label{theo} Let $q=1/p$, with $p\in
\mathbb{N}\backslash\{0,1\}$; let
$\rho=10\pi^2/(5\pi^2-24)$.  Then $\zeta_{q}(2)$
is irrational, and the inequality
\[\left|\zeta_{q}(2)-\frac{a}{b}\right|\leq \left|b\right|^{-\rho}\]
has at most finitely many integer solutions $\left(a,b\right)$.
\end{theorem}

Another way of putting this statement is by using the
\textit{irrationality measure} $\mu$.  There are a number of
equivalent definitions for this irrationality measure
(Liouville-Roth constant, order of approximation, irrationality
exponent). One of them is:
\[\mu\left(x\right)=\inf\{\rho:\left|x-\frac{a}{b}\right|>\frac{1}{b^{\rho+\epsilon}},
 \forall\epsilon>0,\forall a,b \in\mathbb{Z}, b \:\rm{sufficiently}\: \rm{large}\}.\]
Notice that for a rational number $x$ we have $\mu(x)=1$, whereas
for an irrational number $x$ we have $\mu(x)\geq 2$ (\cite{hardy},
Theorem 187). So the theorem implies that
\[  2\leq\mu\left(\zeta_{q}(2)\right)\leq
\frac{10\pi^2}{5\pi^2-24}\approx 3.8936.   \]
This is sharper than the upper bound $4.07869374$ given by Zudilin \cite{zud}.

The proof of the theorem is a $q$-adaptation of proofs of the
irrationality of $\zeta(2)$, as given by Ap\'{e}ry
\cite{Apery}, based on Hermite-Pad\'e approximation \cite{Beukers, wva}.
It uses type I Hermite-Pad\'{e} approximation to two
functions $f_1$ and $f_2$, with the property that
$f_1(1)=\zeta_{q}(1)$ and
$f_2(1)=\zeta_{q}(2)$. The polynomials that
arise in this approximation can be found explicitly because they
are closely related to a specific family of orthogonal polynomials,
namely the little $q$-Jacobi polynomials.

We can prove the irrationality and give the upper bound for the
measure of irrationality of $\zeta_{q}(2)$ using the
following two lemmas:
\begin{lemma} \label{irr}
Let $x$ be a real number and suppose there exist
integer sequences $a_n,b_n,\ (n\in \mathbb{N})$ such that
\begin{enumerate}
    \item $b_nx-a_n\neq 0$ for all $n\in \mathbb{N}$;
    \item $\lim_{n\rightarrow\infty}\left(b_nx-a_n\right)=0.$
\end{enumerate}
Then $x$ is irrational.
\end{lemma}
\begin{lemma}\label{irrmaat}
If the conditions of the previous lemma hold, with
$\left|b_nx-a_n\right|=\mathcal{O}(1/b_n^s)$ and
$b_n<b_{n+1}<b_n^{1+o\left(1\right)}$, then $\mu(x)\leq
1+1/s$.
\end{lemma}
For the latter, see, e.g., \cite[Ex.~3, p.~376]{bor2}.

\begin{remark}\label{oeps}
\rm Since an easy calculation shows that
for a natural number $r$,
\[\sum_{k=1}^\infty \frac{kq^{k}}{1-q^k}-\sum_{k=1}^\infty \frac{kq^{rk}}{1-q^k}
=\sum_{i=1}^{r-1}\frac{q^i}{(1-q^i)^2}\in\mathbb{Q},\]
we also obtain the irrationality of the series
\[   \sum_{k=1}^{\infty}\frac{kq^{rk}}{1-q^k}.\]
Moreover, these numbers obviously have the same irrationality measure as
$\zeta_q(2)$.
\end{remark}

\section{Some $q$-calculus}
The following elements of $q$-calculus will often be used :
\begin{itemize}
    \item The $q$-Pochhammer symbols
    \begin{math}\left(a;q\right)_n=\prod_{j=0}^{n-1}\left(1-aq^j\right)\end{math}
and
    \begin{math}\left(a;q\right)_\infty=\prod_{j=0}^{\infty}\left(1-aq^j\right)\end{math};
    \item The $q$-binomial factors:
\[   \qbinom{n}{k}_q=\frac{\left(q;q\right)_n}{\left(q;q\right)_k\left(q;q\right)_{n-k}} \]
 and
\[  \qbinom{n}{k}_p=\frac{\left(p;p\right)_n}{\left(p;p\right)_k\left(p;p\right)_{n-k}}, \]
where a straightforward calculation shows that
\[   \qbinom{n}{k}_q=q^{k(n-k)} \qbinom{n}{k}_p; \]
    \item The $q$-derivative
\[ D_qf\left(x\right)=\begin{cases}
 \displaystyle \frac{f\left(x\right)-f\left(qx\right)}{x(1-q)}, & \textrm{if }  x\neq 0, \\
                f'(0), & \textrm{if } x=0;
                      \end{cases}  \]
    \item The $q$-integral
\[ \int_0^{q^i}f\left(x\right)d_qx=\sum_{k=i}^\infty
    q^kf(q^{k}) \]
and
\[ \int_{q^j}^{q^i}f(x)\,d_qx=\int_{0}^{q^i}f(x)\,d_qx-\int_{0}^{q^j}f(x)\,d_qx; \]
    \item The $q$-Leibniz rule:
\[  D_q^n[f(x)g(x)]=\sum_{k=0}^n \qbinom{n}{k}_q D_q^k\left(f(x)\right)D_q^{n-k}\left(g(q^kx)\right). \]
\end{itemize}

In the literature, the $q$-integral is often defined with an extra
factor $1-q$, which makes the $q$-integration and the
$q$-derivation inverse operations.  Since we do not need this
property, we drop the factor to prevent it from arising everywhere
in the approximations and the analysis. We will need the little
$q$-Jacobi polynomials. These are given by the explicit formula
(see \cite{Koekoek})
\begin{equation}
 \label{qjacexp}
P_n(x;a,b|q)={}_2\phi_1\left( \left.\begin{array}{c}q^{-n},abq^{n+1} \\ aq \end{array} \right| q;qx \right)
  = \sum_{k=0}^n
  \frac{(q^{-n};q)_k (abq^{n+1};q)_k}{(q;q)_k(aq;q)_k}q^kx^k\ ,
\end{equation}
and there also exists a Rodrigues formula, given by
\begin{equation} \label{qjacrod}
\frac{(qx;q)_\infty}
{(q^{\beta+1}x;q)_\infty}x^\alpha
P_n\left(x;q^\alpha,q^\beta|q\right)
=\frac{q^{n\alpha+n(n-1)/2}(1-q)^n}{(q^{\alpha+1};q)_n}
D_p^n\left(\frac{(qx;q)_\infty}{(q^{\beta+n+1}x;q)_\infty}x^{\alpha+n}\right).
\end{equation}
The orthogonality is given by
\begin{equation}\label{orthojacobi}
\sum_{k=0}^\infty\frac{(bq;q)_k}{(q;q)_k}(aq)^kP_n\left(q^k;a,b|q\right)q^{km}
= 0, \qquad m=0,...,n-1.
\end{equation}

For $q$-integrals there exists an analogue to integration by
parts, which is called summation by parts.

\begin{lemma}\label{sbp}\textsc{[Summation by parts]}
If $g(p)=0$ of $f(1)=0$ and if both series converge, then
\[  \sum_{k=0}^\infty q^k \left. f(q^k)D_pg(x)\right|_{q^k}
 =-q\sum_{k=0}^\infty q^k \left. g(q^k)D_qf(x)\right|_{q^k}.\]
\end{lemma}

\section{Approximations to $\zeta_q(2)$}
\subsection{The Hermite-Pad\'{e} approximation problem}

The following Hermite-Pad\'{e} approximation problem is
considered: find polynomials $A_n$, $B_n$ of degree $\leq n$ and $C_n$
of degree $\leq n-1$, for which
\begin{eqnarray}
   \label{fnpl} F_n(z)=A_n(z)+B_n(z)\log_{q}(z) &=& 0, \qquad \textrm{for } z=1,p,p^2,...,p^n, \\
   \label{pade} A_n(z)f_1(z)+B_n(z)f_2(z)-C_n(z)&=& \mathcal{O}\left(\frac{1}{z^{n+1}}\right),
  \qquad z\rightarrow\infty,
\end{eqnarray}%
with $\log_q z = \log z/\log q$ the logarithm with base $q$ and
\[   f_1(z)=\int_0^1\frac{d_{q}x}{z-x}, \quad
     f_2(z)=\int_0^1\log_{q}x\frac{d_{q}x}{z-x}. \]

We suggest the following expression for $F_n(x)$, where
$R_n$ is a yet unknown polynomial of degree $n$ and $x$ is any
point on the grid $\{q^i |\ i\in\mathbb{Z}\}$:
\begin{equation}\label{fnxint}
  F_n(x) =
\int_x^1R_n\left(\frac{x}{t}\right) (qt;q)_n\frac{d_{q}t}{t}.
\end{equation}
This choice of $F_n$ satisfies the first condition
(\ref{fnpl}): if we use the definition of the $q$-integral we get
\begin{eqnarray}
  \nonumber F_n(p^\ell) &=& -\sum_{k=-\ell}^{-1}
  R_n(p^{\ell+k})(q^{k+1};q)_n   \\
  &=& -\sum_{k=0}^{\ell-1}
  R_n(p^{k})(q^{k+1-\ell};q)_n.
\end{eqnarray}
It is now obvious that for $1\leq \ell\leq n$ we have
$(q^{k+1-\ell};q)_n=(1-q^{k+1-\ell}) \cdots (1-q^{k+n-\ell})$,
and since the summation index $k$ runs from $0$ to $\ell-1$, the
factor $(1-q^0)$ is present in every term, hence
$F_n(p^\ell)=0$.  If $\ell=0$, then $F_n(p^\ell)$ is an empty sum, hence also zero.

The orthogonality conditions that follow from the second part of
the Hermite-Pad\'{e} approximation problem, are given by
\begin{equation}\label{orthoconds}
\int_0^1 F_n(x)x^m\,d_{q}x =0, \qquad m=0,...,n-1,
\end{equation}
and they will allow us to find what the polynomials $R_n$ really
are:
\begin{eqnarray}
 \nonumber 0 &=& \int_0^1 F_n(x)x^m\, d_{q}x \\
 \nonumber  &=& \int_0^1\int_x^1 R_n\left(\frac{x}{t}\right)(qt;q)_n\frac{d_{q}t}{t}x^m\ d_{q}x \\
 \label{nul}  &=& q^{m+1}\int_0^1(qt;q)_nt^m\, d_{q}t \
 \int_0^1R_n(qy)y^m\,d_{q}y.
\end{eqnarray}

\begin{remark}
\rm The diligent reader may think that there
is a mistake in this last expression: by switching the order
of integration and substituting $x$ by $yt$, one would expect a
factor $R_n(y)$ in the last line.  However, if one replaces the
integrals by sums according to the definition of $q$-integration,
and one then switches the order of summation, this factor turns
out to be $R_n(qy)$.  This shows that one has to be extremely
careful when working with $q$-integrals.
\end{remark}

The first integral of (\ref{nul}) is the integral over $[0,1]$ of
a strictly positive integrand, so this factor is positive, which
means that the other integral has to be zero. Comparing this to
(\ref{orthojacobi}), we conclude that the polynomials $R_n$
should really be the little
$q$-Jacobi polynomials $P_n\left(px;1,1|q\right)$ (up to a constant factor),
which are in fact little $q$-Legendre polynomials \cite[\S 3.12.1]{Koekoek}.

We use the explicit expression for the little $q$-Jacobi
polynomials (\ref{qjacexp}) and insert it in (\ref{fnxint}), to
get
\begin{eqnarray}
 \nonumber F_n(x) &=& \int_x^1(qt;q)_n
  \sum_{k=0}^n\frac{(q^{-n};q)_k(q^{n+1};q)_k}{(q;q)_k^2}\frac{x^k}{t^k}\frac{d_{q}t}{t} \\
  \nonumber &=& \sum_{k=0}^n (-1)^kp^{\frac{k(k+1)}{2}-nk}
  \qbinom{n}{k}_{p}\qbinom{n+k}{k}_{p} x^k \int_x^1  \frac{(qt;q)_n}{t^{k+1}}\,d_{q}t.
\end{eqnarray}
Since
\[ (qt;q)_n=\sum_{k=0}^np^{\frac{k(k-1)}{2}-nk}\qbinom{n}{k}_{p}(-1)^k t^k,\]
we only need an expression for
$\int_x^1t^m\, d_{q}t$ with $m=-n-1,...,n-1$.
It is easily shown that
\[  \begin{cases}
   \displaystyle \int_x^1t^m,d_qt=\frac{1-x^{m+1}}{1-q^{m+1}}, & \textrm{if $m\neq-1$;} \\[12pt]
  \displaystyle  \int_x^1t^md_qt=\log_qx, & \textrm{if $m=-1$.}
\end{cases}
\]
Hence
\begin{eqnarray*}
    F_n(x)&=&\sum_{k=0}^n(-1)^kp^{\frac{k(k+1)}{2}-nk}\qbinom{n}{k}_{p}\qbinom{n+k}{k}_{p}
 \sum_{i=0,i\neq k}^n \qbinom{n}{i}_{p}p^{\frac{i(i+1)}{2}-ni}(-1)^i\frac{x^k-x^i}{p^i-p^k} \\
    & &+\ \sum_{k=0}^np^{-2kn+k^2} \qbinom{n}{k}_{p}^2 \qbinom{n+k}{k}_{p}x^k\log_{q}x.
\end{eqnarray*}
So, using the definition of $F_n$ in (\ref{fnpl}), we find the
polynomials $A_n$ and $B_n$ to be
\begin{equation}\label{anx}
A_n(x)=\sum_{k=0}^n(-1)^kp^{\frac{k(k+1)}{2}-nk}\qbinom{n}{k}_{p}\qbinom{n+k}{k}_{p}
 \sum_{i=0,i\neq k}^n \qbinom{n}{i}_{p}p^{\frac{i(i+1)}{2}-ni}(-1)^i\frac{x^k-x^i}{p^i-p^k},
\end{equation}
\begin{equation}\label{bnx}
B_n(x)=\sum_{k=0}^np^{-2kn+k^2}\qbinom{n}{k}_{p}^2 \qbinom{n+k}{k}_{p}x^k.
\end{equation}
The Hermite-Pad\'{e} approximation theory also gives us an
expression for the third polynomial $C_n$.  We have
\begin{equation}\label{cnpn2}
C_n(x)=\sum_{\ell=0}^\infty q^\ell\frac{A_n(x)-A_n(q^\ell)}{x-q^\ell}
+\sum_{\ell=0}^\infty \ell q^\ell\frac{B_n(x)-B_n(q^\ell)}{x-q^\ell}.
\end{equation}
Plugging in the formulae (\ref{anx})--(\ref{bnx}) for $A_n$ and $B_n$  and changing the order of summation,
we arrive at
\begin{eqnarray}
\nonumber C_n(x)&=&\sum_{k=0}^n\sum_{i=0,i\neq k}^n(-1)^{k+i}
\qbinom{n}{k}_{p}\qbinom{n}{i}_{p} \qbinom{n+k}{k}_{p} \frac{p^{-nk-ni+\frac{k(k+1)}{2}+\frac{i(i+1)}{2}}}{p^i-p^k} \\
\nonumber & &\times \sum_{\ell=0}^\infty q^\ell\frac{q^{\ell i}-q^{\ell k}+x^{k}-x^{i}}{x-q^\ell} \\
\nonumber & &+\ \sum_{k=0}^n\qbinom{n}{k}_{p}^2 \qbinom{n+k}{k}_{p}p^{-2kn+k^2}
 \sum_{\ell=0}^\infty\frac{\ell q^\ell}{x-q^\ell}\left(x^{k}-q^{\ell k}\right).
\end{eqnarray}
Using three times the identity
\[\frac{A^n-B^n}{A-B}=\sum_{t=0}^{n-1}A^tB^{n-t-1}\]
we can further isolate the infinite sums, which can now be
calculated.  To that end we use the series
\[\sum_{\ell=0}^\infty \ell A^\ell=\frac{A}{(1-A)^2}, \qquad |A| < 1  \]
in the second term.  This leads us to a closed formula for $C_n$:
\begin{eqnarray}
\nonumber C_n(x)&=&\sum_{k=0}^n\sum_{i=0,i\neq k}^n(-1)^{k+i}
\qbinom{n}{k}_{p} \qbinom{n}{i}_{p} \qbinom{n+k}{k}_{p} \frac{p^{-nk-ni+\frac{k(k+1)}{2}+\frac{i(i+1)}{2}}}{p^i-p^k}\\
\nonumber & &\times \left[\sum_{t=0}^{k-1} \frac{p^{k-t}x^t}{p^{k-t}-1}
 -\sum_{t=0}^{i-1} \frac{p^{i-t}x^t}{p^{i-t}-1}\right]
\\
\label{cnpneindig} & &+\ \sum_{k=0}^n \qbinom{n}{k}_{p}^2 \qbinom{n+k}{k}_{p} p^{-2kn+k^2}
\sum_{t=0}^{k-1} \frac{p^{k-t}x^t}{(p^{k-t}-1)^2}.
\end{eqnarray}
Evaluating (\ref{pade}) at $p^n$ and using (\ref{fnpl}) shows that
${a_n^*}/{b_n^*}$ is an approximation for
$\zeta_q(2)$, with
\begin{eqnarray}\label{anster}
a_n^*&=&B_n(p^n)\sum_{k=1}^{n-1}\frac{k}{p^k-1}+C_n(p^n), \\
\label{bnster} b_n^*&=&B_n(p^n).
\end{eqnarray}

\subsection{Integer sequences}

To get some results regarding irrationality and the irrationality
measure, the numerator and denominator of the approximant should
be integers. So we will have to multiply them by an expression
$e_n$, in such a way that the numbers
\begin{equation}\label{bnan}
a_n=e_na_n^*\qquad{\rm and}\qquad b_n=e_nb_n^*
\end{equation}
are integers.
Looking at the explicit formulae (\ref{bnx})--(\ref{cnpneindig}) for $B_n(p^n)$ and
$C_n(p^n)$, we can deduce the factors that are needed in $e_n$.  Keep in mind that
$p=1/q$ is a natural number greater
than 1.

It is well-known that the $p$-binomial factors are integers, hence
only powers of $p$ can arise in the denominator of $B_n(p^n)$.
There is a factor $p^{k^2-nk}$ in $B_n(p^n)$, with the summation
index $k$ going from $0$ to $n$.  The minimum of this exponent is
obviously $-\left\lfloor \frac{n^2}{4}\right\rfloor$ (with
$\lfloor \cdot \rfloor$ the floor function), so we conclude that
\[p^{\lfloor\frac{n^2}{4}\rfloor}B_n(p^n)\in\mathbb{Z}.\]
The possible denominators that arise in the term
\[B_n(p^n)\sum_{k=1}^{n-1}\frac{k}{p^k-1}\]
are then clearly cancelled out by $p^{\lfloor
\frac{n^2}{4}\rfloor}\ {\rm lcm}\{p^j-1,1\leq j\leq n-1\}$, where lcm denotes the least common multiple.

Finally, we need to find the denominator of $C_n(p^n)$.
 This denominator consists of factors $p^j$ and $p^j-1$.  Looking at (\ref{cnpneindig}), we
 see that the latter are completely cancelled by the factor
\[\left({\rm lcm}\{p^j-1,1\leq j\leq n\}\right)^2.\]
It is well-known that, as a polynomial in $x$,
\[{\rm lcm}\{x^j-1,1\leq j\leq n\}=d_n(x),\]
with
\begin{equation}\label{dnp}
d_n(x)=\prod_{d=1}^n\Phi_d\left(x\right),
\end{equation}
where $\Phi_d$ are the cyclotomic polynomials
defined by
\[ \Phi_d(x)=\prod_{\substack{k=1\\ \gcd(k,d)=1}}^d (x-\omega_d^k), \]
with
\[\omega_d=e^{2\pi i/d}.\]
Hence, putting a factor $d_n^2(p)$ in $e_n$ will cancel all
factors of type $p^j-1$.  The only other denominators that can originate from $C_n(p^n)$, are powers of $p$.

At a first glance, the factor needed to cancel these powers of
$p$, would be $p^{n^2-n}$.  However, calculations using
\texttt{Maple} indicate that a factor $p^{\lfloor \frac{n^2}{4}\rfloor}$ is enough.
This will indeed be proved in the next section.
This leads us to the following
\begin{lemma}\label{enconj}
If we choose
\begin{equation}\label{en}
e_n =p^{\lfloor\frac{n^2}{4}\rfloor}d_{n}^2(p),
\end{equation}
then $a_n$ and $b_n$, as defined in (\ref{bnan}), are integers.
\end{lemma}

\begin{remark}
\rm From (\ref{cnpn2}) or (\ref{cnpneindig}), we see that $C_n(p^n)$ consists
of two terms: one originating from $A_n$, the other originating
from $B_n$.  It is easy to show that both these terms have the
predicted denominator $p^{n^2-n}$, but as mentioned before,
putting these two terms together results in the disappearance of
these high power denominators, making $\left\lfloor \frac{n^2}{4}
\right\rfloor$ the largest remaining exponent of $p$ in the
denominator of $C_n(p^n)$.
\end{remark}

\subsection{Proof of Lemma~\ref{enconj}}

We will work with the $q$-Mellin transform of the expression
$F_n$.  The $q$-Mellin transform of a measurable function $f$
on the $q$-exponential lattice on $(0,1]$ is given by
\[\hat{f}(s)=\int_0^1f(x)x^s\,d_qx.\]
The particular structure of $F_n$ as given in (\ref{fnpl}), and
the orthogonality conditions as stated in (\ref{orthoconds}),
allow us to give an explicit expression for $\widehat{F}_n$:
\[\widehat{F}_n(s)=\frac{(p;p)_n}{p^{n^2+n+1}}\frac{q^s(q^{s-n+1};q)_n}{(q^{s+1};q)_{n+1}^2}.\]
The Hermite-Pad\'{e} theory gives us an expression for the error
term of the approximation.  In this case we get
\[b_n^*\zeta_q(2)-a_n^*=\sum_{k=0}^\infty\frac{q^k}{p^n-q^k}F_n(q^k)
=q^n\sum_{\ell=0}^\infty q^{n\ell}\widehat{F}_n(\ell).\]
Let us now introduce the rational function
\[R_n(T;q)=\frac{T^n(Tq^{-n+1};q)_n}{(qT;q)_{n+1}^2}\]
and the series
\[S_n(q)=\sum_{\ell=0}^\infty q^\ell R_n(q^\ell;q).\]
Then obviously
\begin{equation}\label{snq}
S_n(q)=\frac{p^{n^2+n+1}}{(p;p)_n}\sum_{\ell=0}^\infty
q^{n\ell}\widehat{F}_n(\ell)
=\frac{p^{n^2+2n+1}}{(p;p)_n}\left(b_n^*\zeta_q(2)-a_n^*\right).
\end{equation}
A partial fraction decomposition gives
\[R_n(T;q)=\sum_{s=1}^2\sum_{j=1}^{n+1}\frac{d_{s,j,n}(q)}{\left(1-q^jT\right)^s}\]
with
\[d_{s,j,n}(q)=(-1)^sq^{js}\left.\frac{d^{2-s}}{dT^{2-s}}\left(R_n(T;q)\left(T-q^{-j}\right)^2\right)\right|_{T=q^{-j}}\]
for $s=1,2$.  Isolation of the infinite sums allows us to
recognize the expressions for $\zeta_q(1)$ and $\zeta_q(2)$, and
we obtain
\[S_n(q)=\sum_{j=1}^{n+1}d_{1,j,n}(q)q^{-j}\zeta_q(1)
+\sum_{j=1}^{n+1}d_{2,j,n}(q)q^{-j}\zeta_q(2)-D_1(n,q)-D_2(n,q)\]
with
\[D_1(n,q)=\sum_{j=1}^{n+1}d_{1,j,n}(q)q^{-j}\sum_{\ell=1}^{j-1}\frac{q^\ell}{1-q^\ell} \]
and
\[  D_2(n,q)=\sum_{j=1}^{n+1}d_{2,j,n}(q)q^{-j}\sum_{\ell=1}^{j-1}\frac{q^\ell}{(1-q^\ell)^2}.\]
Since we already know from (\ref{snq}) that $S_n(q)$ is a $\mathbb{Q}$-linear
combination of 1 and $\zeta_q(2)$, and since the three numbers 1,
$\zeta_q(1), \zeta_q(2)$ are $\mathbb{Q}$-linearly independent
(see \cite{kelly}), we see that the coefficient of $\zeta_q(1)$
has to be zero.  Moreover, looking at $D_1$ and $D_2$ and replacing
$q$ by $1/p$, it is an easy task to see that both these quantities
contain a power $p^{\left\lceil\frac{3n^2}{4}\right\rceil+2n+1}$
as a factor in their numerator.  Together with (\ref{snq}) this
allows us to conclude that the highest possible exponent of $p$ in
the denominator of $a_n^*$ is
$\left\lfloor\frac{n^2}{4}\right\rfloor$, and hence that the
factor $e_n$ as proposed in (\ref{en}) indeed makes $a_n$ and
$b_n$ integers.
\medskip

Since we have an explicit expression for the coefficient of
$\zeta_q(1)$, its vanishing yields a $q$-binomial identity:
\begin{corollary}
The following identity holds
\[\sum_{j=0}^nq^{-2nj+j^2}\qbinom{n+j}{n}_q \qbinom{n}{j}_q^2
\left(n+\sum_{k=1}^{n+j}\frac{1}{1-q^k}-3\sum_{k=1}^j\frac{1}{1-q^k}+2\sum_{k=1}^{n-j}\frac{q^k}{1-q^k}\right) = 0.\]
Multiplying by $1-q$ and letting $q$ tend to 1, we obtain the identity:
\[\sum_{j=0}^n\binom{n+j}{n}\binom{n}{j}^2\left(H(n+j)+2H(n-j)-3H(j)\right) = 0. \]
where $H(n)=\sum_{k=1}^n 1/k$ are harmonic numbers.
\end{corollary}

\section{Irrationality of \textbf{$\zeta_{q}(2)$}}
\subsection{Estimate for the error term}

So far we know that $a_n$ and $b_n$ are integers, and that
${a_n}/{b_n}$ is an approximation of
$\zeta_{q}(2)$. Now we want to estimate
$\left|b_n\zeta_{q}(2)-a_n\right|$.  To meet the
conditions of Lemma \ref{irr}, we need to prove that this quantity
tends to zero as $n$ tends to infinity, and that it is never
zero.
Once again we use the expression for the error term of the
approximation:
\begin{equation} \label{bnzq2minan}
b_n\zeta_{q}\left(2\right)-a_n
=e_n\sum_{k=0}^\infty\frac{q^k}{p^{n}-q^k}F_n(q^k).
\end{equation}
In this last expression we need $F_n(q^k)$.  This can
be calculated using (\ref{fnxint}):
\begin{equation}\label{fouttermalsint}
F_n(q^k)=\int_{q^k}^1R_n\left(\frac{q^k}{t}\right)(qt;q)_n\ \frac{d_{q}t}{t}
=\sum_{\ell=0}^{k-1}P_n\left(q^{k-\ell-1};1,1|q\right)(q^{\ell+1};q)_n.
\end{equation}
If we now use the Rodrigues formula for the little $q$-Jacobi
polynomials (\ref{qjacrod}) and plug this into
(\ref{fouttermalsint}), then after changing the order of summation we
get
\begin{equation}\nonumber
\left|b_n\zeta_{q}(2)-a_n\right|=e_n\frac{q^{\frac{n(n-1)}{2}+1}(1-q)^n}{(q;q)_n}
\left|\sum_{\ell=0}^\infty (q^{\ell+1};q)_nq^\ell
\sum_{k=0}^\infty \frac{q^k}{p^n-q^{k+\ell+1}}D_{p}^n\left.\left[(qx;q)_nx^{n}\right]\right|_{x=q^k}\right|.
\end{equation}
Applying $n$ times summation by parts (Lemma \ref{sbp}) we have
\begin{eqnarray*}
\left|b_n\zeta_{q}(2)-a_n\right|&=&e_n\frac{q^{\frac{n(n+1)}{2}+1}(1-q)^n}{(q;q)_n}\\
& &\times \left|\sum_{\ell=0}^\infty (q^{\ell+1};q)_nq^\ell
\sum_{k=0}^\infty (q^{k+1};q)_nq^kq^{nk}D_{q}^n\left.\left(\frac{1}{p^n-q^{\ell+1}x}\right)\right|_{x=q^k}\right|.
\end{eqnarray*}
Now it can be proven by induction that the $q$-derivative needed
in this last expression, is given by
\[D_{q}^n\frac1{p^n-q^{\ell+1}x}
=\frac{q^{\ell n}(q;q)_np^{\frac{n(n-1)}{2}}}{(1-q)^n\prod_{j=0}^n\left(p^{n+j}-xq^{\ell+1}\right)}.\]
We recognize a double $q$-integral for
$\left|b_n\zeta_{q}(2)-a_n\right|$:
\begin{equation}\nonumber
\left|b_n\zeta_{q}(2)-a_n\right|=
e_nq^{n+1} \left|\int_0^1\int_0^1\frac{(qx;q)_nx^n(qy;q)_ny^n}{\prod_{j=0}^n (p^{n+j}-qxy)}\,d_{q}x\,d_{q}y\right|.
\end{equation}
None of these factors is zero, and the integrand is strictly
positive on $(0,1]^2$, so we see that the first
condition of Lemma \ref{irr} is fulfilled.

Obviously $\prod_{j=0}^n(p^{n+j}-qxy)$
reaches its minimum in $(x,y)=(1,1)$.
Moreover, the function $(qx;q)_nx^n$ reaches its
maximum in $x=1$, as long as $0<q\leq\frac12$, which is the case
we are working with since $p=1/q$ is an integer. To see this,
it is enough to show that $x(1-q^mx)/(1-q^m)\leq 1$ for all $x\in[0,1]$
and for $m=1,...,n$. So we can make the
estimate
\begin{eqnarray}
 \nonumber \left|b_n\zeta_{q}(2)-a_n\right| &\leq& e_n\frac{(q;q)_n^2}{(1-q)^2}
\frac{q^{n+1}}{\prod_{j=0}^n(p^{n+j}-q)} \\
 \label{errorfinaal}  &=&e_n\ \frac{(q;q)_n^2q^{n+1}}{(1-q)^2(q^{n+1};q)_{n+1}}q^{\frac{3}{2}n(n+1)}.
\end{eqnarray}

\subsection{Asymptotic behaviour}

The asymptotic behaviour of the cyclotomic polynomials is known
(see e.g. \cite{lqlp}) and is given in the following lemma.
\begin{lemma}\label{dnasymp}
Suppose $p$ is an integer greater than one and let
$d_n$ be given by (\ref{dnp}).  Then
\[\lim_{n\rightarrow\infty}d_n(p)^{1/n^2}=p^{3/\pi^2}. \]
\end{lemma}

Hence the expression (\ref{en}) has the asymptotic behaviour
\begin{equation}\label{enasymp}
\lim_{n\rightarrow\infty} e_n^{1/n^2}=p^{\frac{6}{\pi^2}+\frac14}
\end{equation}
and (\ref{errorfinaal}) has the asymptotic behaviour
\begin{equation} \label{bnalspmacht}
\lim_{n\rightarrow\infty}\left|b_n\zeta_{q}(2)-a_n\right|^{1/n^2}\leq
p^{\frac{6}{\pi^2}+\frac14-\frac32}\approx p^{-0.6421}.
\end{equation}
So we conclude that also the second
condition of Lemma \ref{irr} is fulfilled, and hence that
$\zeta_{q}(2)$ is irrational.

\begin{remark}
\rm One could try to use the same method to prove the irrationality of
\[ \zeta_{q_1,q_2}(2)=\sum_{k=1}^\infty \frac{kq_1^k}{1-q_2^k}\]
with $q_2=1/p_2, q_1=1/p_1$ and integers $p_1,p_2$.  Little $q$-Jacobi
polynomials with different parameters are needed in this case.
However, the $e_n$ which is needed to cancel the denominators,
turns out to be too large and
\[\lim_{n\rightarrow\infty}|b_n\zeta_{q_1,q_2}(2)-a_n|^{1/n^2}>1.\]
Hence we cannot deduce the irrationality for this family of
numbers.  The case where $p_1$ and $p_2$ are related in a certain
way (they are both powers of the same integer $p$) gives
asymptotically better results, but still not good enough to prove
irrationality.  So we only obtain an irrationality result for the
family of numbers mentioned in Theorem~\ref{theo} and
Remark~\ref{oeps}.
\end{remark}

\subsection{The measure of irrationality}

To use Lemma \ref{irrmaat}, we need to get a value for $s$ in
$\left|b_n\zeta_{q}(2)-a_n\right|=\mathcal{O}\left(1/{b_n^s}\right)$.
We already know that
\[\lim_{n\rightarrow\infty} \left|b_n\zeta_{q}(2)-a_n\right|^{1/n^2}\leq
p^{\frac{6}{\pi^2}-\frac{5}{4}}.\]
If we can now find the
asymptotic relation between $b_n=e_nB_n(p^n)$ and
$p^{n^2}$, then we obtain the desired value for $s$.  From the
explicit formula (\ref{bnx}) for $B_n(p^n)$ and the
asymptotic behaviour of $e_n$ in (\ref{enasymp}), it is clear that
\[\lim_{n\rightarrow\infty}b_n^{1/n^2}=p^{\frac{6}{\pi^2}+\frac14+1}=p^{\frac{24+5\pi^2}{4\pi^2}}.\]
Together with (\ref{bnalspmacht}), this means that
\[\left|b_n\zeta_{q}(2)-a_n\right|=\mathcal{O}\left(1/{b_n^\frac{5\pi^2-24}{5\pi^2+24}}\right).\]
Hence Lemma \ref{irrmaat} gives us an upper bound for the
measure of irrationality:
\[\mu\left(\zeta_{q}(2)\right)\leq1+\frac{5\pi^2+24}{5\pi^2-24}=\frac{10\pi^2}{5\pi^2-24},\]
which concludes the proof of Theorem \ref{theo}.  Moreover,
\texttt{Maple} calculations indicate that the sequence
$\log_{b_n}\left(|b_n\zeta_q(2)-a_n|\right)$ indeed takes values
that are close to $-\frac{5\pi^2-24}{5\pi^2+24}$.\\

\begin{verbatim}
Department of Mathematics
Katholieke Universiteit Leuven
Celestijnenlaan 200B - box 2400
BE-3001 Leuven
Belgium
christophe@wis.kuleuven.be
walter@wis.kuleuven.be
\end{verbatim}

\end{document}